\documentclass{amsart}
\pdfoutput=1

\usepackage{amsfonts, amssymb, amsmath, amsthm}                               
\usepackage{slashed}
\usepackage{graphicx}
\usepackage{caption}
\usepackage{subcaption}

\newtheorem{theorem}{Theorem} 	      	      	                              
\newtheorem{corollary}[theorem]{Corollary}     	      	      	      	      
\newtheorem{proposition}[theorem]{Proposition} 	      	      	      	      
\newtheorem{definition}[theorem]{Definition} 	      	      	                
\newtheorem*{remark}{Remark}                                                  

\numberwithin{equation}{section}                                              
\numberwithin{theorem}{section}                                               

\newcommand{\ol}[1]{\overline{#1}}                                            
\newcommand{\mc}[1]{\mathcal{#1}}                                             

\newcommand{\R}{\mathbb{R}}                                                   

\newcommand{\grad}{\nabla^\sharp}                                             
\newcommand{\nasla}{\slashed{\nabla}}                                         

\newcommand{\beq}{\begin{equation}}
\newcommand{\eeq}{\end{equation}}

\def\XXint#1#2#3{{\setbox0=\hbox{$#1{#2#3}{\int}$ }
\vcenter{\hbox{$#2#3$ }}\kern-.6\wd0}}

\allowdisplaybreaks

\begin{document}


\title[Subconformal focusing nonlinear waves]{On the profile of energy concentration at blow-up points for subconformal focusing nonlinear waves.}

\author{Spyros Alexakis}
\address{Department of Mathematics\\ 
University of Toronto\\
40 St George Street Rm 6290
\\ Toronto, ON M5S 2E4\\ Canada}
\email{alexakis@math.utoronto.ca}

\author{Arick Shao}
\address{Department of Mathematics\\
South Kensington Campus\\
Imperial College\\
London SW7 2AZ\\ United Kingdom}
\email{c.shao@imperial.ac.uk}

\begin{abstract}
We consider singularities of the focusing subconformal nonlinear wave equation and some generalizations of it.
At noncharacteristic points on the singularity surface, Merle and Zaag have identified the rate of blow-up of the $H^1$-norm of the solution inside cones that terminate at the singularity.
We derive bounds that restrict how this $H^1$-energy can be distributed inside such cones.
Our proof relies on new localized estimates---obtained using Carleman-type inequalities---for such nonlinear waves. These  bound the $L^{p+1}$-norm  in the interior of timelike cones by their $H^1$-norm near the boundary of the cones.
Such estimates can also be applied to obtain certain integrated decay estimates for globally regular solutions to such equations, in the interior of time cones.
\end{abstract}

\maketitle



\section{Introduction} \label{sec:intro}

The main aim of this paper is to study the behaviour near a singularity of  solutions to focusing power-law-type wave equations of the form
\begin{align}
\label{focNLW} \Box \phi + V (t,x) |\phi|^{p-1} \phi = 0 \text{,} \qquad \Box := - \partial_t^2 + \Delta_x \text{,}
\end{align}
over the Minkowski space-time $\mathbb{R}^{n+1}$,
\begin{align*}
(t,x) \in \R \times \R^n \text{,} \qquad V (t, x) \geq \epsilon > 0 \text{,}
\end{align*}
with $p$ lying in the \emph{subconformal} range,
\begin{align*}
p \in \left( 1, 1 + \frac{4}{n-1} \right) \text{.}
\end{align*}

A detailed understanding of this has been obtained in one spatial dimension; see \cite{CF, CZ, MZ1Da, MZ1Db, MZ1Dc, MZ1Dd}.
Our focus here is when $n \ge 2$, without symmetry assumptions.
  
Our starting point are the works of Merle and Zaag, \cite{MZ1, MZ2, MZ3}, who obtained (sharp) upper and lower bounds on the blow-up rate of the $H^1$-norm of the solution inside cones that terminate at the singularity. 
These results, however do not give information on the concentration properties of the $H^1$-norm along the different directions toward the singularity.
We show a result of this type, proving that any fraction of the $H^1$-energy that persists on certain specific balls that shrink down towards the singularity must be ``distributed'' over such balls. 
Slightly more precisely, we show that at any space-time blow-up point $O$ for \eqref{focNLW} and any time-like cone $\mc{C}$ that terminates at $O$, the fraction of the $H^1$-norm that lives inside $\mc{C}$ cannot concentrate in a thinner cone $\mc{C}'\subset\mc{C}$.
In fact, we prove that the amount of $H^1$-energy in $\mc{C} \setminus \mc{C}'$ bounds the energy in $\mc{C}'$. 
\medskip  
 
To put our result in context, we recall a few known facts from the literature.
One can easily construct blow-up solutions for the classical focusing nonlinear wave equations---\eqref{focNLW} with $V \equiv 1$---by imposing spatial homogeneity of the solution.
This results in an ODE in $t$, for which one finds that the function
\begin{align*}
\phi_* (t, x) := C (-t)^\frac{-2}{p-1} \text{,} \qquad C := \left[ \frac{ 2 (p+1) }{ (p-1)^2 } \right]^\frac{1}{p-1}
\end{align*}
solves \eqref{focNLW} and blows up (for all $x \in \mathbb{R}^n$) at $t=0$.
While this function clearly fails to be in $L^2$, one can produce examples of singularity formation starting from $L^2$-initial data by applying a spatial cut-off function to $\phi_* (t, x)$ at $t = -1$ for $|x|\ge M$, $M > 0$ large.
By the domain of dependence property of \eqref{focNLW}, we then still obtain a singularity at $t = 0, x = 0$ (and nearby).

One can then build more blow-up solutions based on this example by using the invariance of the equation under Lorentz transformations. 
This yields the family of the solutions $\kappa(d,y)$.\footnote{$d$ encodes the Lorentz transformations that fixes the blow-up point; $y$ stands for the similarity variable, which we do not use here.}
A rich literature has been developed in understanding the universality of these solutions.

It is known---see \cite{Alinhac, MZ1}---that if a solution of \eqref{focNLW} blows up at one point, then it must blow up on a Lipschitz hypersurface $\Gamma$ that is graph over $\{ t = 0 \}$, with
\begin{align*}
\Gamma := (T(x), x) \text{,} \qquad x \in \mathbb{R}^n \text{,} \qquad | T(x) - T(y) | \leq | x - y | \text{.}
\end{align*}
This is called the \emph{blow-up graph} of the solution. 
 
\begin{figure} \label{non-char}
\centering
\begin{subfigure}[b]{0.4\textwidth}
\includegraphics[width=130pt]{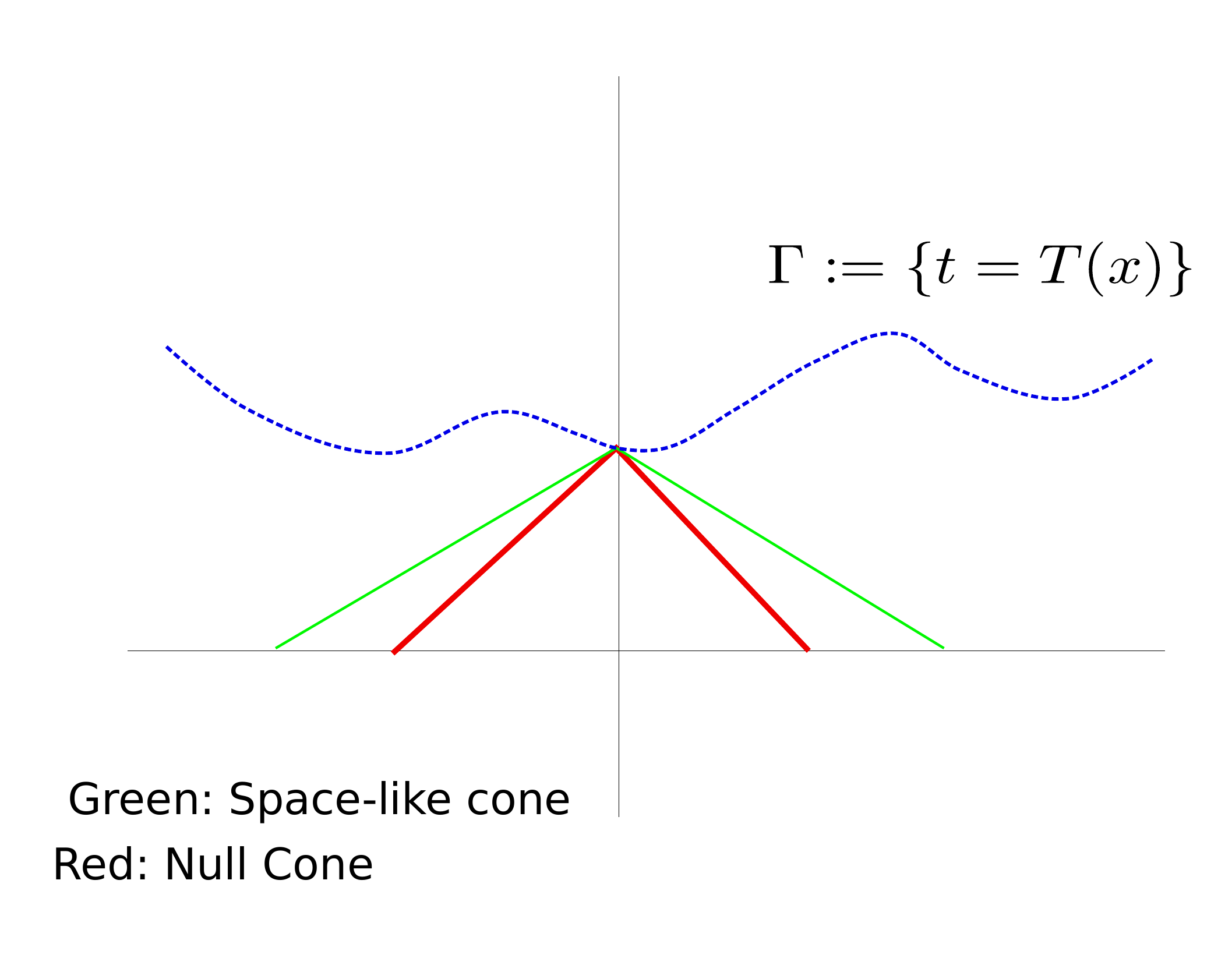}   
\caption{The blow-up graph $\Gamma$ and a non-characteristic point on $\Gamma$.}
\label{fig.non-char}
\end{subfigure}
\nobreak
\qquad\qquad\qquad
\begin{subfigure}[b]{0.3\textwidth}
\includegraphics[height=100pt]{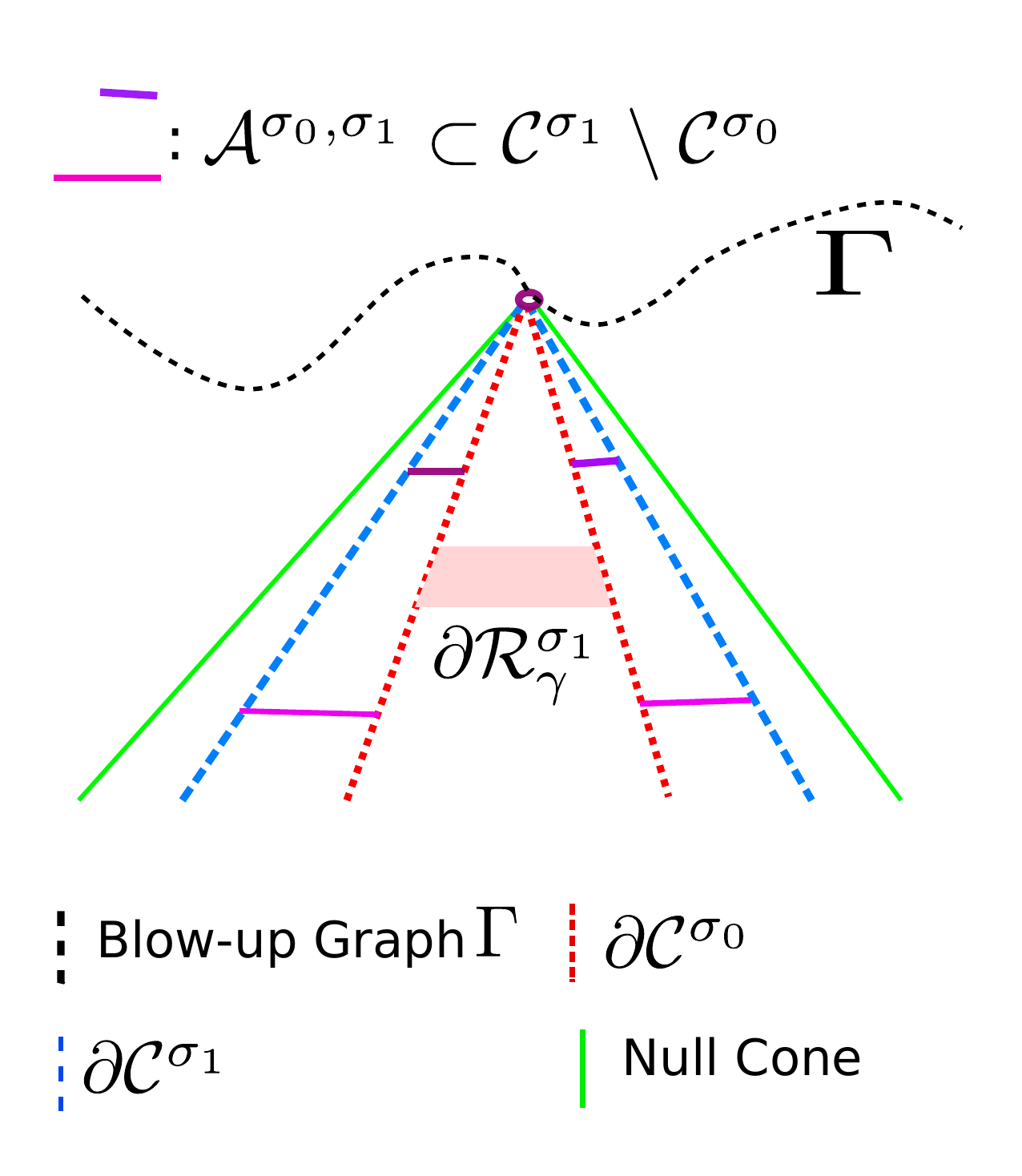}
\label{fig.annuli}
\caption{The annuli $\mc{A}^{\sigma_0,\sigma_1}(t)$ in $\mc{C}^{\sigma_1}\setminus\mc{C}^{\sigma_0}$}
\end{subfigure}
\end{figure}

Understanding the profile of the blow-up at each point of $\Gamma$ involves a distinction between \emph{characteristic} and \emph{noncharacteristic} points.
The latter (in the language adopted in our paper) are points $P$ for which there exists a backwards spacelike cone\footnote{More precisely, $\mc{C}$ is invariant under dilations centered at $P$, and the restriction of the Minkowski metric onto $\mc{C}$ is spacelike.} $\mc{C}$ emanating from $P$ so that $\mc{C}$ locally intersects $\Gamma$ only at $P$; see figure \ref{fig.non-char}.
In particular, any point of $\Gamma$ that differentiable and spacelike is noncharacteristic.

At noncharacteristic points on $\Gamma$, Merle and Zaag, \cite{MZ1}, showed (in the case $V \equiv 1$) that the rate of blow-up of the $H^1$-norm when approaching from a backwards null cone matches the one provided by the ODE examples $\phi_\ast$.
Translating so that the singularity is at the origin, they prove that if $(0, 0) \in \Gamma$ is noncharacteristic, then there exists an $\varepsilon > 0$ so that
\begin{align}
\label{MZ} \varepsilon &\le (-t)^\frac{2}{p-1} \frac{ \| \phi (t) \|_{ L^2 ( B (0, - t) ) } }{ (-t)^\frac{n}{2} } + (-t)^{\frac{2}{p-1} + 1} \frac{ \| \partial_t \phi (t) \|_{ L^2 ( B (0, -t) ) } }{ (-t)^\frac{n}{2} } \\
\notag &\qquad + (-t)^{\frac{2}{p-1} + 1} \frac{ \| \nabla_x \phi (t) \|_{ L^2 ( B(0, - t  ) ) } }{ (-t)^\frac{n}{2} } \text{.}
\end{align}
Also, for any $\sigma \in (0,1)$ (even if $(0, 0) \in \Gamma$ is characteristic), there exists $K_\sigma$ so that
\begin{align}
\label{MZ'} &(-t)^\frac{2}{p-1} \frac{ \| \phi (t) \|_{ L^2 ( B (0, - \sigma t) ) } }{ (-t)^\frac{n}{2} } + (-t)^{\frac{2}{p-1} + 1} \frac{ \| \partial_t \phi (t) \|_{ L^2 ( B (0, -\sigma t) ) } }{ (-t)^\frac{n}{2} } \\
\notag &\qquad + (-t)^{\frac{2}{p-1} + 1} \frac{ \| \nabla_x \phi (t) \|_{ L^2 ( B(0, -\sigma t  ) ) } }{ (-t)^\frac{n}{2} } \leq K_\sigma \text{.}
\end{align}

\begin{remark}
In fact, the lower bound \eqref{MZ} also holds near noncharacteristic points for more general $V$, as long as $V$ can be extended beyond the singularity.
This can be seen via a straightforward adaptation of \cite[Lemma 3.1]{MZ1}.
In particular, even with these more general $V$, the local well-posedness theory for \eqref{focNLW} in $H^1$ remains the same as for $V \equiv 1$; see, for example, \cite{LS}.
\end{remark}
    
We note that very recent results, \cite{donn_schor:blowup_subcritical, MZ2, MZ3}, also show the stability of the family of blow-up profiles $\kappa (d, y)$ (in similarity variables), as well as the $\mc{C}^1$-regularity of the blow-up graph near such a blow-up profile, for noncharacteristic points.
We also note the results of Killip, Stovall, and Vi\c{s}an, \cite{kill_stov_vis:blowup_nkg}, which apply to nonlinear wave and Klein-Gordon equations in the full energy-subcritical range, $p \in (1, 1 + \frac{4}{n - 2})$; among other things, they established bounds on the local energy inside a null cone emanating from a point on the blow-up surface.

Nonetheless, to our knowledge, there does not seem to be a general result describing the profile of the singularities in higher dimensions beyond these perturbations of the ODE examples $\phi_\ast$.
In fact, Biz\'on, Breitenlohner, Maison, and Wasserman, \cite{BBMZ}, demonstrate in certain cases that there exists a countable family of self-similar solutions, which shows the richness of potential profiles.

In particular, \eqref{MZ} and \eqref{MZ'} do not give information on the profile of the solution in the balls $B(0, -\sigma t )$, but, merely on the growth of its $H^1$ norm. 
For example, the bound \eqref{MZ} in principle allows for the solution $\phi$ to blow up slowly on the annular regions $B (0, -\sigma t) \setminus B (0, -\sigma' t)$ for some $\sigma' < \sigma$, and for the bulk of the $H^1$-energy to be concentrated in the smaller balls $B(0, -\sigma' t)$.
We show that this cannot occur.

To explain our result, we first introduce some notation:

\begin{definition} \label{gradient}
For a $\mc{C}^1$-function $\phi$ defined on some domain in $\R^{n+1}$, we let $\nabla \phi$ denote the spacetime gradient of $\phi$.
In addition, we let $| \nabla \phi |$ denote its vector norm:
\begin{align*}
| \nabla \phi |^2 := ( \partial_t \phi )^2 + \sum_{i = 1}^n ( \partial_{ x^i } \phi )^2 \text{.}
\end{align*}
\end{definition}

\begin{definition} \label{sound.cone}
Given $0 < \sigma < 1$, we let $\mc{C}^\sigma$ the interior of the future time cone about the time axis and with vertex at the origin:
\begin{align*}
\mc{C}^\sigma := \{ (t, x) \mid 0 < r < \sigma t \} \text{,} \qquad r := |x| \text{.}
\end{align*}
\end{definition}

\begin{definition} \label{annuli}
Given $0 < \sigma_0 < \sigma_1 < 1$ and $t \in \R$, we define the annulus
\begin{align*}
\mc{A}^{\sigma_0, \sigma_1} (t) := \{ (t, x) \mid \sigma_1 |t| < r < \sigma_2 |t| \} \text{.}
\end{align*}
\end{definition}

\begin{definition}
\label{R}
Given any $t_*\ne 0$ and any $\gamma > 1$, we let
\begin{align}
\mc{R}_\gamma^\sigma (t_*) := \begin{cases} \{ \gamma t_* < t < \gamma^{-1} t_* \} \cap ( - \mc{C}^\sigma ) & t_\ast < 0 \text{,} \\ \{ \gamma^{-1} t_* < t < \gamma t_* \} \cap \mc{C}^\sigma & t_\ast > 0 \text{.} \end{cases}
\end{align}
\end{definition} 

Note that for the ODE examples $\phi_*$ (for $V \equiv 1$), we have:
\begin{align}
\label{bl-up1} (-t)^{-n + 2 + \frac{4}{p-1}} \int_{ \mc{A}^{\sigma_0,\sigma_1} (t) } [ (\partial_t \phi_*)^2 + (\partial_r\phi_*)^2 ] &= C_{\sigma_0, \sigma_1} \text{,} \\
\notag (-t)^{-n + \frac{4}{p-1}} \int_{ \mc{A}^{\sigma_0, \sigma_1} (t) } |\phi_*|^2 &= C'_{\sigma_0, \sigma_1} \text{,}
\end{align}
Similarly, on the domains $\mc{R}^{\sigma_1}_\gamma(t)\subset\mc{C}_\sigma$, we see:
\begin{align}
\label{bl-up2} (-t)^{-n + 1 + \frac{4}{p-1}} \int_{ \mc{R}^{\sigma_0}_\gamma (t) } [ (\partial_t \phi_*)^2 + (\partial_r\phi_*)^2 ] &= C_{\sigma_0, \gamma} \text{,} \\
\notag (-t)^{-n - 1 + \frac{4}{p-1}} \int_{ \mc{R}^{\sigma_0}_\gamma (t) } |\phi_*|^2 &= C'_{\sigma_0, \gamma} \text{.}
\end{align}

We show that if any fraction of the (weighted) $H^1$-norm in $\mc{C}$ moves in timelike directions towards the singularity, then the norm must be ``smeared'' across the various angles that terminate at the singularity:

\begin{theorem}
\label{singularity.foc}
Let $0 < \sigma_0 < \sigma_1 < 1$.
Assume $\phi$ is a $\mc{C}^2$-solution on the domain
\begin{align}
\label{eq.singularity_ass_C} \mc{C} := \{ (t, x) \in \R^{n + 1} \mid -1 < t < 0 \text{, } |x| < \sigma_1 |t| \} \text{.}
\end{align}
of \eqref{focNLW}, where $p$ and $V$ satisfy
\begin{align}
\label{eq.singularity_ass_pV} 1 \leq p < 1 + \frac{4}{n - 1} \text{,} \qquad V \in \mc{C}^1 ( \bar{\mc{C}} ) \text{,} \qquad V \gtrsim 1 \text{.}
\end{align}
Suppose also that
\begin{align}
\label{eq.singularity_ass} 0 \leq \limsup_{t \rightarrow 0^-} |t|^{2 - n + \frac{4}{p-1}} \int_{ \mc{A}^{\sigma_0, \sigma_1} (t) } ( | \nabla \phi |^2 + |t|^{-2} |\phi|^2 ) < \delta \text{.}
\end{align}
Then, there is some $\gamma > 1$ (also depending on $\sigma_0$ and $\sigma_1$) such that
\begin{align}
\label{eq.singularity} \limsup_{t \rightarrow 0^-} |t|^{1 - n + \frac{4}{p-1} } \int_{ \mc{R}^{\sigma_0}_\gamma (t) } ( | \nabla \phi |^2 + |t|^{-2} \phi^2 ) \lesssim \delta \text{.}
\end{align}
\end{theorem}

A corollary of this is that any fraction of the (weighted) $H^1$-norm of $\phi$ in the balls $B(0, -t)$ in \eqref{MZ} cannot be concentrated entirely on smaller balls $B(0, - \sigma t)$:
 
\begin{corollary} \label{cor.singularity.foc}
Consider any solution $\phi$ of \eqref{focNLW} on the domain $\mc{C}$, defined as in \eqref{eq.singularity_ass_C}, where $p$ and $V$ satisfy \eqref{eq.singularity_ass_pV}.
Assume also that $(0,0)$ belongs to the singular set of $\phi$, and assume there exists $0 < \sigma_0 < 1$ so that 
\begin{align}
\label{some.in} \limsup_{ t \to 0^-} |t|^{2 - n + \frac{4}{p-1}} \int_{ |x| < \sigma_0 |t| } [ |\nabla \phi|^2 + |t|^{-2} \phi^2 ] > 0 \text{.}
\end{align}
Then, for every $\sigma_0 < \sigma_1 < 1$,
\begin{align}
\label{side-blowup} \limsup_{t \to 0^-} \left[ |t|^{2 - n + \frac{4}{p-1}} \int_{ \mc{A}^{\sigma_0, \sigma_1} (t)} |\nabla \phi|^2 + |t|^{-n + \frac{4}{p-1}} \int_{ \mc{A}^{\sigma_0, \sigma_1} (t) } |\phi|^2 \right] > 0 \text{.}
\end{align}
\end{corollary}
 
We remark that the analysis in the papers \cite{MZ1, MZ2, MZ3, MZ1Da, MZ1Db, MZ1Dc, MZ1Dd} relies crucially on the existence of a Lyapunov functional for the equation in similarity variables.
This is especially useful in completely understanding the blow-up behaviour in one 
spatial dimension.
Our methods here do not rely on any Hamiltonian structure or on a Lyapunov functional.
We also do not work directly in similarity variables, although our approach of 
studying time cones that terminate at the singularity perhaps captures the same 
geometric information as the similarity variables.
For this reason, the results we derive directly apply to the more general class of equations \eqref{focNLW}, with a non-constant $V(t,x)$.
We remark that the class of equations to which our results apply can be further broadened, but we do not pursue this for brevity.
\footnote{In particular, the Carleman estimates of \cite{alex_shao:uc_global} can be generalized to other wave operators; see \cite[Sect. 2]{alex_shao:uc_global} for further discussions.}

Our method relies on a new estimate derived from the Carleman-type estimates in \cite{alex_shao:uc_global}; for focusing equations the latter precisely works in the subconformal range of the exponents $p$.
This new estimate is stated precisely in Section \ref{sec:intro_est} and is proven in Section \ref{sec:est}.
Note the subconformal and conformal range of exponents is precisely that for which the energy functional $E$ used in \cite{MZ3} (in the similarity variables) is well-defined and monotone decreasing.
 
\subsection{Integrated decay in time-cones.}

We present a second result, which deals with decay in time of solutions of \eqref{focNLW} inside time cones.
We show that integrated decay \emph{on} the cone $\partial \mc{C}^\sigma$, captured by the finiteness of the $H^1$-norm, implies integrated decay \emph{inside} the cone, measured by the finiteness of the bulk $L^{p+1}$-norm.

We note that for technical reasons, the result
in this subsection does not apply to all subconformal exponents, but rather to those for which 
\begin{align}
\label{range} p \in \left( 1, 1 + \frac{4n}{n^2+n-4} \right) \text{,} \qquad n \geq 2 \text{.}
\end{align}
This agrees with the full subconformal range only when $n=2$.
In higher dimensions, this is a strictly smaller range of exponents.
We remark, however, that one can obtain the same 
results in the full subconformal range 
in all dimensions, by strengthening the decay assumption on the cone,
measured in a suitable $L^q$-norm, with $q$ depending on $n$.

\begin{theorem} \label{out-in-decay1}
Let $0 < \sigma < 1$.
Assume $\phi$ is a $\mc{C}^2$-solution on
\begin{align*}
\bar{\mc{C}}^{\sigma, \ast} := \ol{ \mc{C}^\sigma \cap \{ t > 1 \} }
\end{align*}
of \eqref{focNLW}, where $p$ satisfies,
\begin{align}
\label{eq.final_state_ass_p} 1 \leq p < \begin{cases} 1 + \frac{4}{n - 1} & n \leq 2 \text{,} \\ 1 + \frac{4n}{n^2 + n - 4} & n \geq 3 \text{.} \end{cases}
\end{align}
and where $V \in \mc{C}^1 ( \bar{\mc{C}}^{\sigma, \ast} )$ satisfies
\begin{align}
\label{eq.final_state_ass_V} V \simeq 1 \text{,} \qquad | \nabla V | = o ( t^{-1} ) \text{,}
\end{align}
the latter as $t \nearrow \infty$.
If $\phi$ satisfies
\begin{align} \label{eq.final_state_ass} \int_{ \partial \mc{C}^\sigma \cap \{ t > 1 \} } ( | \nabla \phi |^2 + | \phi |^{p + 1} ) < \infty \text{,}
\end{align}
then
\begin{align}
\label{eq.final_state} \int_{ \mc{C}^\sigma \cap \{ t > 1 \} } t^{-1} |\phi|^{p+1} < \infty \text{.}
\end{align}
\end{theorem}

\begin{remark}
Recall that any time cone in $\R^{n+1}$ (around any point and time-axis) can be mapped onto some $\mc{C}^\sigma$ via a Minkowski isometry.
Consequently, an analogue of Theorem \ref{out-in-decay1} can be readily derived for any time cone.
\end{remark}

\subsection{Localized Estimates in Time Cones} \label{sec:intro_est}

The proofs of Theorems \ref{singularity.foc} and \ref{out-in-decay1} rely on new estimates for wave equations on suitable finite domains inside time cones.
We feature two such estimates here, as we believe these may be of independent interest.
In general, these time-localized estimates are proved using the global Carleman estimates developed in \cite{alex_shao:uc_global}.

\begin{definition} \label{cones}
Let $\sigma > 0$ and $\eta > 1$.
Then, for any $t_\ast > 0$, we define
\begin{align*}
\mc{K}^\sigma_\eta (t_\ast) := \partial \mc{C}^\sigma \cap \{ \eta^{-1} t_\ast < t < \eta t_\ast \} \text{.}
\end{align*}
\end{definition}


Both estimates control a solution $\phi$ of \eqref{focNLW} in a time slab $\mc{R}^\sigma_\gamma (t_\ast)$ inside the cone.
In the first estimate, this is controlled by $\phi$ and $\nabla \phi$ in a slab $\mc{K}^\sigma_\eta (t_\ast)$ of the bounding time cone, while in the second estimate, this bound is by its values in a family of annuli $\mc{A}^{\sigma, \sigma'} (\tau)$ for $\tau \simeq t_*$.
The precise estimates are below:
 
\begin{theorem} \label{thm.main_est_timecone}
Fix $0 < \sigma < 1$ and $t_\ast > 0$, and suppose $\eta > 1$ is sufficiently large.
In addition, assume $\phi \in \mc{C}^2 ( \ol{ \mc{R}^\sigma_\eta (t_\ast) } )$ is a solution of \eqref{focNLW}, where
\begin{align}
\label{eq.main_est_timecone_ass_p} 1 \leq p < \begin{cases} 1 + \frac{4}{n - 1} & n \leq 2 \text{,} \\ 1 + \frac{4n}{n^2 + n - 4} & n \geq 3 \text{,} \end{cases}
\end{align}
and where $V \in \mc{C}^1 ( \ol{ \mc{R}^\sigma_\eta (t_\ast) } )$ satisfies
\begin{align}
\label{eq.main_est_ass_V} C^{-1} \leq V \leq C \text{,} \qquad | \nabla V | \leq \alpha t_\ast^{-1} \text{,}
\end{align}
for some $C > 1$ and some $\alpha > 0$, sufficiently small with respect to $p$ and $C$.
Then, there exists $1 < \gamma < \eta$, with $\gamma - 1$ sufficiently small, such that
\begin{align}
\label{eq.main_est_timecone} \int_{ \mc{R}^\sigma_\gamma (t_*) } |\phi|^{p+1} \lesssim t_\ast \int_{ \mc{K}^\sigma_\eta (t_*) } [ | \nabla \phi |^2 + |\phi|^{p+1} + t_\ast^{-2} | \phi |^2 ] \text{.}
\end{align}
\end{theorem}
  
\begin{theorem} \label{thm.main_est_annulus}
Fix $0 < \sigma_0 < \sigma_1 < 1$ and $t_\ast > 0$, and suppose $\eta > 1$ is sufficiently large.
In addition, assume $\phi \in \mc{C}^2 ( \mc{R}^{\sigma_1}_\eta (t_\ast) )$ is a solution of \eqref{focNLW}, where
\begin{align}
\label{eq.main_est_annulus_ass_p} 1 \leq p < 1 + \frac{4}{n - 1} \text{.}
\end{align}
and where $V \in \mc{C}^1 ( \mc{R}^{\sigma_1}_\eta (t_\ast) )$ satisfies \eqref{eq.main_est_ass_V} for some $C > 1$ and some $\alpha > 0$, small with respect to $p$ and $C$.
Then, there is $1 < \gamma < \eta$, with $\gamma - 1$ small, such that
\begin{align}
\label{eq.main_est_annulus} \int_{ \mc{R}^{\sigma_0}_\gamma (t_*) } |\phi|^{p+1} \lesssim t_\ast \sup_{ \eta^{-1} t_\ast \leq \tau \leq \eta t_\ast } \int_{ \mc{A}^{\sigma_0, \sigma_1} (\tau) } ( | \nabla \phi |^2 + |\phi|^{p+1} + t_\ast^{-2} | \phi |^2 ) \text{.}
\end{align}
\end{theorem}

We will apply both estimates in this paper: Theorem \ref{thm.main_est_annulus} is used to prove Theorem \ref{singularity.foc}, while Theorem \ref{thm.main_est_timecone} is used for Theorem \ref{out-in-decay1}.
However, this choice is not essential; one could prove a variant of Theorem \ref{singularity.foc} using Theorem \ref{thm.main_est_timecone} (and altering the assumptions accordingly), and the same statement holds for Theorems \ref{out-in-decay1} and \ref{thm.main_est_annulus}.
For conciseness, we only demonstrate one application for each estimate.



Note that Theorem \ref{thm.main_est_timecone} has the advantage that one bounds by the values of $\phi$ on only a single hypersurface.
On the other hand, Theorem \ref{thm.main_est_annulus} has the advantange that it applies to all subconformal $p$ in all dimensions.
In particular, an analogue of Theorem \ref{out-in-decay1} using \eqref{eq.main_est_annulus} would apply for all subconformal $p$.

\section{The Main estimates} \label{sec:est}

In this section, we establish the main estimates, Theorems \ref{thm.main_est_timecone} and \ref{thm.main_est_annulus}.

\subsection{The Global Carleman Estimate}

The starting point of the proof is the (nonlinear) global Carleman estimate of \cite[Theorem 2.18]{alex_shao:uc_global}, in the focusing ($+$) case.
For the reader's convenience, we restate this estimate, in this particular case, below.
For this, we first need a few preliminary definitions from \cite[Sect. 2]{alex_shao:uc_global}:

\begin{definition} \label{def.D}
Let $\mc{D}$ denote the exterior of the null cone about the origin,
\begin{align*}
\mc{D} := \{ Q \in \R^{n+1} \mid r (Q) > |t| (Q) \} \text{,}
\end{align*}
and let $f$ denote the Lorentz square distance from the origin:
\begin{align*}
f \in \mc{C}^\infty (\mc{D}) \text{,} \qquad f := \frac{1}{4} ( r^2 - t^2 ) \text{.}
\end{align*}
\end{definition}

Note in particular that $f$ is strictly positive on $\mc{D}$.
Next, we define regions within $\mc{D}$ for which the Carleman estimate applies.
These are precisely the regions for which one can apply the standard divergence theorem in the Lorentzian setting.

\begin{definition} \label{def.region_adm}
We say that an open, connected subset $\Omega \subseteq \mc{D}$ is \emph{admissible} iff:
\begin{itemize}
\item The closure of $\Omega$ is a compact subset of $\mc{D}$.

\item The boundary $\partial \Omega$ of $\Omega$ is piecewise smooth, with each smooth piece being either a spacelike or a timelike hypersurface of $\mc{D}$.
\end{itemize}
For an admissible $\Omega \subseteq \mc{D}$, we define the \emph{oriented unit normal} $\mc{N}$ of $\partial \Omega$ as follows:
\begin{itemize}
\item $\mc{N}$ is the inward-pointing unit normal on each spacelike piece of $\partial \Omega$.

\item $\mc{N}$ is the outward-pointing unit normal on each timelike piece of $\partial \Omega$.
\end{itemize}
Integrals over such an admissible region $\Omega$ and portions of its boundary $\partial \Omega$ will be with respect to the volume forms induced by $g$.
\end{definition}

We can now state the global Carleman estimate, \cite[Theorem 2.18]{alex_shao:uc_global}:

\begin{theorem} \label{thm.carleman_nl}
Let $\phi \in \mc{C}^2 (\bar{\Omega})$, with $\Omega \subseteq \mc{D}$ be an admissible region.
Furthermore, let $p \geq 1$, and let $V \in \mc{C}^1 (\bar{\Omega})$ be strictly positive.
\footnote{While $\phi$ and $V$ were assumed to be regular on $\mc{D}$ in \cite{alex_shao:uc_global}, the entire proof takes place on $\bar{\Omega}$.}
Then, for any $a > 0$,
\begin{align}
\label{eq.carleman_nl} \frac{1}{p + 1} \int_\Omega f^{2a} \cdot V \Gamma_V \cdot | \phi |^{p + 1} &\leq \frac{1}{8 a} \int_\Omega f^{2 a} f \cdot | \Box_V \phi |^2 + \int_{ \partial \Omega } P^V_\beta \mc{N}^\beta \text{,}
\end{align}
where $\mc{N}$ is the oriented unit normal to $\partial \Omega$, and where:
\begin{align}
\label{eq.bulk_gamma} \Box_V \phi &:= \Box \phi + V | \phi |^{p - 1} \phi \text{,} \\
\notag \Gamma_V &:= \nabla^\alpha f \nabla_\alpha (\log V) - \frac{n - 1 + 4 a}{4} \left( p - 1 - \frac{4}{n - 1 + 4a} \right) \text{,} \\
\notag P^V_\beta &:= f^{2 a} \left( \nabla^\alpha f \cdot \nabla_\alpha \phi \nabla_\beta \phi - \frac{1}{2} \nabla_\beta f \cdot \nabla^\mu \phi \nabla_\mu \phi \right) \\
\notag &\qquad + \frac{1}{p + 1} f^{2 a} \nabla_\beta f \cdot V | \phi |^{p + 1} + \left( \frac{n - 1}{4} + a \right) f^{2 a} \cdot \phi \nabla_\beta \phi \\
\notag &\qquad + a \left( \frac{n - 1}{4} + a \right) f^{2 a} f^{-1} \nabla_\beta f \cdot \phi^2 \text{.}
\end{align}
\end{theorem}

\begin{remark}
One particularly useful feature of \eqref{eq.carleman_nl} is that the weight $f^{2a}$ vanishes as one approaches the double null cone about the origin.
\end{remark}

\subsection{A Preliminary Estimate}

Consider now the interior $\mc{C}^\sigma$, $0 < \sigma < 1$, of a time cone (see Definition \ref{sound.cone}), and suppose $\phi \in \mc{C}^2 (\bar{\mc{C}}^\sigma)$.
In addition:

\begin{definition}
\label{cuts.cones}
Let $\zeta$ be a future timelike geodesic ray beginning from the origin and lying within $\mc{C}^\sigma$.
\footnote{These correspond to all future unit timelike vectors $\zeta$ at the origin that lie in $\mc{C}^\sigma$.}
Furthermore, given $\zeta$ and $t_\ast$ as above:
\begin{itemize}
\item Let $\zeta (t_\ast)$ denote the intersection point of $\zeta$ with $\{ t = t_\ast \}$.

\item Let $\mc{D}^\sigma_{t_\ast, \zeta}$ denote the spacetime region lying both in the exterior of the double null cone emanating from $\zeta (t_\ast)$ and inside $\mc{C}^\sigma$:
\begin{align*}
\mc{D}^\sigma_{t_\ast, \zeta} := \{ Q \in \R^{n+1}\mid | t (Q) - t_\ast | < | x (Q) - x ( \zeta (t_\ast) ) | \text{, } r (Q) < \sigma \cdot t (Q) \} \text{.}
\end{align*}

\item We denote the boundary of $\mc{D}^\sigma_{t_\ast, \zeta}$ on the time cone by
\begin{align*}
\mc{G}^\sigma_{t_\ast, \zeta} := \partial \mc{D}^\sigma_{t_\ast, \zeta} \cap \partial \mc{C}^\sigma \text{.}
\end{align*}

\item We define the translated square distance function
\begin{align*}
f_{t_\ast, \zeta} := \frac{1}{4} [ | x - x ( \zeta (t_\ast) ) |^2 - | t - t_\ast |^2 ] \text{.}
\end{align*}
\end{itemize}
\end{definition}

Our goal here will be to apply Theorem \ref{thm.carleman_nl} to prove a preliminary estimate for $\phi$ in the region $\mc{D}^\sigma_{t_\ast, \zeta}$, for some $\zeta$ and $t_\ast$ as above, by the values of $\phi$ and $\nabla \phi$ on $\mc{G}^\sigma_{t_\ast, \zeta}$.
The result is given in the subsequent proposition:

\begin{proposition} \label{thm.carleman_timecone}
Let $p \in (1, 1 + \frac{4}{n - 1})$, and choose $a > 0$ so that
\begin{align}
\label{eq.carleman_timecone_ass} 1 \leq p < \frac{4}{n - 1 + 4a} \text{,}
\end{align}
Fix also $0 < \sigma < 1$ and $t_\ast > 0$, and assume $V \in \mc{C}^1 ( \ol{ \mc{D}^\sigma_{t_\ast, \zeta} } )$ satisfies the conditions \eqref{eq.main_est_ass_V}, where $C > 1$ and where $\alpha > 0$ is sufficiently small with respect to $p$, $a$, and $C$.
Then, for any solution $\phi \in \mc{C}^2 ( \ol{ \mc{D}^\sigma_{t_\ast, \zeta} } )$ of \eqref{focNLW}, we have the estimate
\footnote{In fact, one can slightly strengthen \eqref{eq.carleman_timecone}: ``$| \nabla \phi |^2$" in the right-hand side can be replaced by ``$( \partial_t \phi )^2 + ( \partial_r \phi )^2$". This can be accomplished via a more careful expansion in \eqref{eql.boundary_timecone_1}.}
\begin{align}
\label{eq.carleman_timecone} \int_{ \mc{D}^\sigma_{t_\ast, \zeta} } f_{t_\ast, \zeta}^{2a} |\phi|^{p+1} &\lesssim t_\ast^{1 + 4a} \int_{ \mc{G}^\sigma_{t_\ast, \zeta} } | \nabla \phi |^2 + t_\ast^{1 + 4a} \int_{ \mc{G}^\sigma_{t_\ast, \zeta} } |
\phi|^{p+1} \\
\notag &\qquad + t_\ast^{-1 + 4a} \int_{ \mc{G}^\sigma_{t_\ast, \zeta} } \phi^{2} + t_* \int_{ \mc{G}^\sigma_{t_\ast, \zeta} } f_{t_\ast, \zeta}^{-1 + 2a} \phi^2  \text{.}
\end{align}
\end{proposition}

The main idea is to apply Theorem \ref{thm.carleman_nl} to the \emph{translates}
\begin{align*}
\phi_{t_\ast, \zeta} := \phi (\cdot + \zeta (t_\ast))
\end{align*}
on the correspondingly shifted region
\begin{align*}
\{ Q - \zeta (t_\ast) \mid Q \in \mc{D}^\sigma_{t_\ast, \zeta} \} \subseteq \mc{D} \text{.}
\end{align*}
By a simple change of variables, this is equivalent to applying Theorem \ref{thm.carleman_nl} to $\phi$ on the region $\mc{D}^\sigma_{t_\ast, \zeta}$, but with each instance of $f$ in \eqref{eq.carleman_nl} and \eqref{eq.bulk_gamma} replaced by $f_{t_\ast, \zeta}$.
In particular, $\grad f$ in \eqref{eq.carleman_nl} and \eqref{eq.bulk_gamma} is replaced by
\begin{align*}
\grad f_{t_\ast, \zeta} = \frac{1}{2} ( t - t_\ast ) \partial_t + \frac{1}{2} \sum_{ i = 1 }^n [ x^i - x^i ( \zeta (t_\ast) ) ] \cdot \partial_{x^i} \text{.}
\end{align*}

\begin{remark}
The essential observation is that the Carleman weight $f^{2a}_{t_\ast, \zeta}$ vanishes on the null portion of the boundary of $\mc{D}^\sigma_{t_\ast, \zeta}$.
The upshot is that in \eqref{eq.carleman_timecone}, the integral of $\phi$ in $\mc{D}^\sigma_{t_\ast, \zeta}$ is controlled only by boundary terms on $\partial \mc{C}^\sigma$, and not within $\mc{C}^\sigma$.
\end{remark}

\subsubsection{Proof of Proposition \ref{thm.carleman_timecone}}

Consider the region
\begin{align*}
\mc{D}^{\sigma, \varepsilon}_{t_\ast, \zeta} := \{ Q \in \mc{D}^\sigma_{t_\ast, \zeta} \mid f_{t_\ast, \zeta} (Q) > \varepsilon \} \text{,}
\end{align*}
for some small, fixed $\varepsilon > 0$.
Applying Theorem \ref{thm.carleman_nl} to $\Omega := \mc{D}^{\sigma, \varepsilon}_{t_\ast, \zeta}$ (and keeping in mind the aforementioned change of variables), we obtain
\begin{align}
\label{eql.carleman_timecone_1} \frac{1}{p + 1} \int_{ \mc{D}^{\sigma, \varepsilon}_{t_\ast, \zeta} } f_{t_\ast, \zeta}^{2a} \cdot V \Gamma_V \cdot | \phi |^{p + 1} &\leq \int_{ \partial \mc{D}^{\sigma, \varepsilon}_{t_\ast, \zeta} } P^V_\beta \mc{N}^\beta \text{,}
\end{align}
where $\mc{N}$ is the outward-pointing unit normal to $\partial \mc{C}_\sigma$.

Since $t \simeq t_\ast$ and $r \lesssim t_\ast$ on $\mc{D}^\sigma_{t_\ast, \zeta}$, we have
\begin{align}
\label{eql.carleman_timecone_20} | \nabla^\alpha f_{t_\ast, \zeta} \nabla_\alpha V | &\lesssim |t - t_\ast| | \partial_t V | + \sum_{i = 1}^n | x^i - x^i ( \zeta (t_\ast) ) | | \partial_{x^i} V | \lesssim t_\ast | \nabla V | \text{.}
\end{align}
Since the right-hand side of \eqref{eql.carleman_timecone_20} is assumed to be small by \eqref{eq.main_est_ass_V} (with respect to $a, p, C$), then the condition \eqref{eq.carleman_timecone_ass} implies that
\begin{align*}
V \Gamma_V &= \nabla^\alpha f_{t_\ast, \zeta} \nabla_\alpha V + \frac{n - 1 + 4 a}{4} \left( 1 + \frac{4}{n - 1 + 4a} - p \right) V \gtrsim 1 \text{,}
\end{align*}
and hence \eqref{eql.carleman_timecone_1} becomes
\begin{align}
\label{eql.carleman_timecone_2} \epsilon_{a, p, C} \int_{ \mc{D}^\sigma_{ t_\ast, \zeta } } f_{t_\ast, \zeta}^{2a} \cdot | \phi |^{p + 1} &\leq \int_{ \partial \mc{D}^\sigma_{ t_\ast, \zeta } } P^V_\beta \mc{N}^\beta \text{.}
\end{align}

It remains to expand and control the boundary integrand $P^V_\beta \mc{N}^\beta$.
Observe that
\begin{align*}
\partial \mc{D}^{\sigma, \varepsilon}_{ t_\ast, \zeta } &= \mc{G}^{\sigma, \varepsilon}_{t_\ast, \zeta} \cup \mc{F}^{\sigma, \varepsilon}_{t_\ast, \zeta} \text{,} \\
\mc{G}^{\sigma, \varepsilon}_{t_\ast, \zeta} &:= \mc{G}^\sigma_{t_\ast, \zeta} \cap \{ f_{t_\ast, \zeta} > \varepsilon \} \text{,} \\
\mc{F}^{\sigma, \varepsilon}_{t_\ast, \zeta} &:= \mc{D}^\sigma_{t_\ast, \zeta} \cap \{ f_{t_\ast, \zeta} = \varepsilon \} \text{.}
\end{align*}
We begin with the portion $\mc{G}^{\sigma, \varepsilon}_{t_\ast, \zeta}$ lying on the time cone.
Note first that
\begin{align*}
\mc{N} = ( 1 - \sigma^2 )^{ -\frac{1}{2} } ( \sigma \partial_t + \partial_r ) \text{.}
\end{align*}
Furthermore, simple computations yield on $\mc{G}^\sigma_{t_\ast, \zeta}$ that
\begin{align} \label{simple.bounds}
| f_{t_\ast, \zeta} | \lesssim t_\ast^2 \text{,} \qquad \mc{N} ( f_{t_\ast, \zeta} ) \simeq \sigma \cdot t_\ast - r ( \zeta (t_\ast) ) \simeq t_\ast \text{.}
\end{align}
Recalling the definition of $P^V$ in \eqref{eq.bulk_gamma} (and replacing $f$ by $f_{t_\ast, \zeta}$), we bound
\begin{align}
\label{eql.boundary_timecone_1} P^V_\beta \mc{N}^\beta &\lesssim f_{t_\ast, \zeta}^{2 a} | \nabla^\alpha f_{t_\ast, \zeta} \nabla_\alpha \phi \cdot \mc{N} \phi | + f_{t_\ast, \zeta}^{2 a} \mc{N} ( f_{t_\ast, \zeta} ) \cdot ( \partial_t \phi )^2 + f_{t_\ast, \zeta}^{2 a} | \phi \cdot \mc{N} \phi | \\
\notag &\qquad + f_{t_\ast, \zeta}^{2 a} | \mc{N} ( f_{t_\ast, \zeta} ) | \cdot | \phi |^{p + 1} + f_{t_\ast, \zeta}^{-1 + 2 a} | \mc{N} ( f_{t_\ast, \zeta} ) | \cdot \phi^2 \\
\notag &\lesssim t_\ast^{1 + 4 a} ( | \nabla \phi |^2 + | \phi |^{p + 1} ) + t_\ast^{4 a} | \phi | ( | \partial_t \phi | + | \partial_r \phi | ) + f_{t_\ast, \zeta}^{-1 + 2 a} t_\ast \cdot \phi^2 \\
\notag &\lesssim t_\ast^{1 + 4 a} ( | \nabla \phi |^2 + | \phi |^{p + 1} ) + t_\ast^{-1 + 4 a} | \phi |^2 + f_{t_\ast, \zeta}^{-1 + 2 a} t_\ast \cdot \phi^2 \text{,}
\end{align}
where we applied \eqref{simple.bounds} in the intermediate steps.

For $\mc{F}^{\sigma, \varepsilon}_{t_\ast, \zeta}$, the relevant observations are
\begin{align}
\label{simpler.bounds} \nabla^\alpha f_{t_\ast, \zeta} \nabla_\alpha f_{t_\ast, \zeta} = f_{t_\ast, \zeta} \text{,} \qquad \mc{N} = - f_{t_\ast, \zeta}^{-\frac{1}{2}} \grad f_{t_\ast, \zeta} \text{.}
\end{align}
Using \eqref{simpler.bounds}, a direct computation yields on $\mc{F}^{\sigma, \varepsilon}_{t_\ast, \zeta}$ that
\footnote{See \cite[Lemma 3.5]{alex_shao:uc_global} for details behind this computation.}
\begin{align}
\label{eq.boundary_nl_est} P^V_\beta \mc{N}^\beta &\lesssim f_{t_\ast, \zeta}^{ 2 a } [ f_{t_\ast, \zeta}^\frac{1}{2} \cdot | \nasla \phi |^2 + (n + a)^2 f_{t_\ast, \zeta}^{-\frac{1}{2}} \cdot \phi^2 ] \text{.}
\end{align}
Since $\phi$ and $\nabla \phi$ are presumed to be uniformly bounded on $\bar{\mc{D}}^\sigma_{t_\ast, \zeta}$, it follows from \eqref{eq.boundary_nl_est}, in addition to \cite[Lemma 3.10]{alex_shao:uc_global}, that
\begin{align}
\label{eql.boundary_timecone_2} \lim_{ \varepsilon \searrow 0 } \int_{ \mc{F}^{\sigma, \varepsilon}_{t_\ast, \zeta} } P^V_\beta \mc{N}^\beta = 0 \text{.}
\end{align}
Finally, combining \eqref{eql.carleman_timecone_2} and \eqref{eql.boundary_timecone_1}, taking the limit $\varepsilon \searrow 0$, and recalling \eqref{eql.boundary_timecone_2} results in the desired estimate \eqref{eq.carleman_timecone}.

\subsection{Proofs of Theorems \ref{thm.main_est_timecone} and \ref{thm.main_est_annulus}}

We are now prepared to prove the main estimates, \eqref{eq.main_est_timecone} and \eqref{eq.main_est_annulus}.
The first steps, discussed below, are identical for both estimates.
Throughout, we assume $t_\ast$, $\eta$, $p$, $V$ as in the hypotheses of Theorems \ref{thm.main_est_timecone} and \ref{thm.main_est_annulus}.
Moreover, in the case of Theorem \ref{thm.main_est_timecone}, we let $\sigma$ be as in the theorem statement, while for Theorem \ref{thm.main_est_timecone}, we let $\sigma \in (\sigma_0, \sigma_1)$.

\begin{figure} 
\includegraphics[width=140pt]{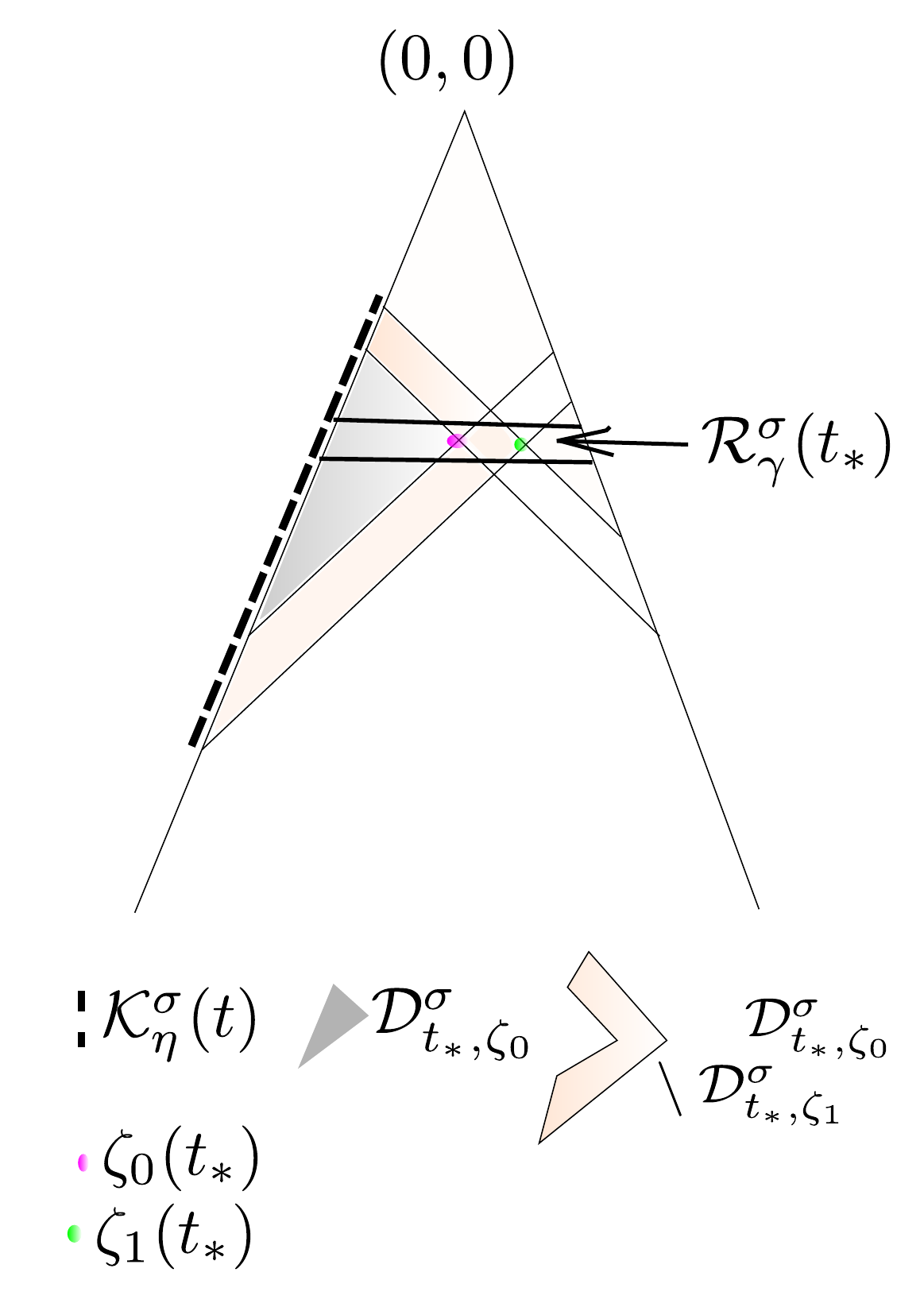}
\caption{Basic setup for the proof of Theorems \ref{thm.main_est_timecone} and \ref{thm.main_est_annulus}.}
\label{overlapping}
\end{figure}

The main step is to apply Proposition \ref{thm.carleman_timecone} twice, with two rays $\zeta_0 \neq \zeta_1$, and to add the resulting estimates.
\footnote{For instance, we can assume $\zeta_0$ is the positive time axis, so that $\zeta_0 (t_\ast) = (t_\ast, 0)$.}
Recalling the slabs $\mc{R}^\sigma_\gamma (t_\ast)$ from Definition \ref{R}, we observe that for some $\gamma > 1$ (depending on $\zeta_0, \zeta_1$),
\begin{align}
\label{eq.domain_1} \mc{R}^\sigma_\gamma (t_*) \subset \mc{D}^\sigma_{t_\ast, \zeta_0} \cup \mc{D}^\sigma_{t_\ast, \zeta_1} \text{.}
\end{align}
(See Figure \ref{overlapping}). 
Moreover, if $\eta > \gamma$ is sufficiently large, then the boundaries of $\mc{D}^\sigma_{t_\ast, \zeta_i}$ are contained in the timelike annuli $\mc{K}^\sigma_\eta (t_*)$ defined in Definition \ref{cones}:
\begin{align}
\label{eq.domain_2} \mc{G}^\sigma_{\zeta_0, t_\ast} \cup \mc{G}^\sigma_{\zeta_1, t_\ast} \subset \mc{K}^\sigma_\eta (t_*) \text{.}
\end{align}
 (Again, see figure \ref{overlapping}). 

Next, for technical reasons, we define $f_{t_\ast, \zeta_i}$ to be zero outside its natural domain of definition, i.e., in the complement of $\mc{D}^\sigma_{t_\ast, \zeta_i}$.
Consider then the function
\begin{align*}
\bar{f}_{t_\ast} := \max (f_{t_\ast, \zeta_0}, f_{t_\ast, \zeta_1}) \in \mc{C} (\R^{n+1}) \text{.}
\end{align*}

We now apply Proposition \ref{thm.carleman_timecone} twice---the first at time $t_\ast$ and direction $\zeta_0$, and the second with $t_\ast$ and $\zeta_1$---with $a \in (0, \frac{1}{2}]$ satisfying \eqref{eq.carleman_timecone_ass}, though its precise value will be determined later.
Adding the two resulting inequalities and recalling the observations \eqref{eq.domain_1} and \eqref{eq.domain_2}, we obtain
\begin{align}
\label{eql.Carl2} \int_{ \mc{R}^\sigma_\gamma (t_\ast) } \bar{f}_{t_\ast}^{2a} |\phi|^{p+1} &\lesssim t_\ast^{1 + 4a} \int_{ \mc{K}^\sigma_\eta (t_\ast) } [ ( \partial_t \phi )^2 + ( \partial_r \phi )^2 + | \phi |^{p + 1} + t_\ast^{-2} | \phi |^2 ] \\
\notag &\qquad + t_\ast \int_{ \mc{G}^\sigma_{t_\ast, \zeta_0} } f_{t_\ast, \zeta_0}^{-1 + 2a} |\phi|^2 + t_\ast \int_{ \mc{G}^\sigma_{t_\ast, \zeta_1} } f_{t_\ast, \zeta_1}^{-1 + 2a} |\phi|^2 \text{.}
\end{align}
Since $\gamma$ is small, $\bar{f}_1 \gtrsim 1$ on $\mc{R}^\sigma_\gamma (1)$, so a rescaling yields $\bar{f}_{t_\ast} \gtrsim t_\ast^2$ on $\mc{R}^\sigma_\gamma (t_\ast)$, hence
\begin{align}
\label{eql.Carl3} t_\ast^{4a} \int_{ \mc{R}^\sigma_\gamma (t_\ast) } |\phi|^{p+1} &\lesssim t_\ast^{1 + 4a} \int_{ \mc{K}^\sigma_\eta (t_\ast) } [ ( \partial_t \phi )^2 + ( \partial_r \phi )^2 + | \phi |^{p + 1} + t_\ast^{-2} | \phi |^2 ] \\
\notag &\qquad + t_\ast \int_{ \mc{G}^\sigma_{t_\ast, \zeta_0} } f_{t_\ast, \zeta_0}^{-1 + 2a} |\phi|^2 + t_\ast \int_{ \mc{G}^\sigma_{t_\ast, \zeta_1} } f_{t_\ast, \zeta_1}^{-1 + 2a} |\phi|^2 \text{.}
\end{align}

Consider the shifted radial and angular parameters,
\begin{align}
\label{eq.radius_angle} r_{t_\ast, \zeta_i} := | x - x ( \zeta_i (t_\ast) ) | \text{,} \qquad \tan \theta_{t_\ast, \zeta_i} := r_{t_\ast, \zeta_i}^{-1} ( t - t_\ast ) \text{.}
\end{align}
Observe that on $\mc{G}^\sigma_{t_\ast, \zeta_i}$, the angular parameter satisfies
\begin{align}
\label{eql.Carl3_rho_0} | \theta_{t_\ast, \zeta_i} | \leq \frac{\pi}{4} \text{.}
\end{align}
Noting that $t \simeq r_{t_\ast, \zeta_i} \simeq t_\ast$ on $\mc{G}^\sigma_{t_\ast, \zeta_i}$, we can estimate
\begin{align}
\label{eql.Carl3_rho_1} f_{t_\ast, \zeta_i}^{-1 + 2a} &\lesssim t_\ast^{-2 + 4 a} ( \tan^{-2} \theta_{t_\ast, \zeta_i} - 1 )^{-1 + 2a} \\
\notag &= t_\ast^{-2 + 4 a} ( \cos 2 \theta_{t_\ast, \zeta_i} )^{-1 + 2 a} ( \sin \theta_{t_\ast, \zeta_i} )^{1 - 2 a} \\
\notag &\lesssim t_\ast^{-2 + 4 a} ( \cos 2 \theta_{t_\ast, \zeta_i} )^{-1 + 2a} \text{,}
\end{align}
where in the last step, we used \eqref{eql.Carl3_rho_0} and that $a \leq 1/2$.

From \eqref{eql.Carl3} and \eqref{eql.Carl3_rho_1}, cancelling factors of $t_\ast^{4a}$, we obtain
\begin{align}
\label{eql.Carl4} \int_{ \mc{R}^\sigma_\gamma (t_\ast) } |\phi|^{p+1} &\lesssim t_\ast \int_{ \mc{K}^\sigma_\eta (t_\ast) } [ ( \partial_t \phi )^2 + ( \partial_r \phi )^2 + | \phi |^{p + 1} + t_\ast^{-2} | \phi |^2 ] \\
\notag &\qquad + t_\ast^{-1} \int_{ \mc{G}^\sigma_{t_\ast, \zeta_0} } ( \cos 2 \theta_{t_\ast, \zeta_0} )^{-1 + 2a} |\phi|^2 \\
\notag &\qquad + t_\ast^{-1} \int_{ \mc{G}^\sigma_{t_\ast, \zeta_1} } ( \cos 2 \theta_{t_\ast, \zeta_1} )^{-1 + 2a} |\phi|^2 \text{.}
\end{align}
From \eqref{eql.Carl4}, we can now derive both \eqref{eq.main_est_timecone} and \eqref{eq.main_est_annulus}.

\subsubsection{Proof of \eqref{eq.main_est_timecone}}

It remains to control the last two terms in the right-hand side of \eqref{eql.Carl4}.
Since both terms are handled in the same manner, we fix $i \in \{ 0, 1 \}$ below, and we abbreviate $\mc{G}^\sigma_{t_\ast, \zeta_i}$ and $\theta_{t_\ast, \zeta_i}$ by $\mc{G}$ and $\theta$, respectively.
\medskip

\noindent
\emph{The case $n \geq 3$:}
Note the assumption \eqref{eq.main_est_timecone_ass_p} for $p$ implies that
\begin{align*}
\frac{1}{2} - \frac{1}{n} < \frac{1}{p - 1} - \frac{n - 1}{4} \text{.}
\end{align*}
Assume $a$ now also satisfies
\begin{align*}
\frac{1}{2} - \frac{1}{n} < a < \min \left( \frac{1}{2}, \frac{1}{p - 1} - \frac{n - 1}{4} \right) \text{.}
\end{align*}
Note this $a$ indeed satisfies \eqref{eq.carleman_timecone_ass}, and moreover,
\begin{align}
\label{eql.ass_a} (-1 + 2a) \frac{n}{2} > -1 \text{.}
\end{align}

An application of H\"older's inequality yields
\begin{align}
\label{eql.Carl4_1} \int_{ \mc{G} } ( \cos 2 \theta )^{-1 + 2a} |\phi|^2 &\lesssim \left[ \int_{ \mc{G} } ( \cos 2 \theta )^{ (-1 + 2a) \frac{n}{2} } \right]^\frac{2}{n} \left( \int_{ \mc{G} } |\phi|^\frac{2n}{n - 2} \right)^\frac{n - 2}{n} \text{.}
\end{align}
By the usual Sobolev embedding $H^1 \hookrightarrow L^\frac{2n}{n-2}$ in $\R^n$, we see that
\begin{align}
\label{eql.Carl4_2} \left( \int_{ \mc{G} } |\phi|^\frac{2n}{n - 2} \right)^\frac{n - 2}{n} &\lesssim \int_{ \mc{G} } ( | \nabla \phi |^2 + t_\ast^{-2} | \phi |^2 ) \text{.}
\end{align}

For the remaining integral, defining $h := \tan \theta$, and applying the coarea formula (i.e., foliating $\mc{G}$ using level sets of $h$), we see that
\begin{align}
\label{eql.Carl4_3} \int_{ \mc{G} } ( \cos 2 \theta )^{ (-1 + 2a) \frac{n}{2} } &\lesssim t_\ast \int_{-1}^1 ( \cos 2 \theta )^{ (-1 + 2a) \frac{n}{2} } \cdot \int_{ \{ h = h' \} } 1 \cdot d h' \text{.}
\end{align}
Note we have on $\mc{G}$ that
\begin{align*}
\frac{d h}{d \theta} = \cos^{-2} \theta \simeq 1 \text{.}
\end{align*}
Moreover, note that each level set of $h$ is a spherical cross-section of $\mc{G}$ and hence has volume $O (t_\ast^{n-1})$.
Consequently, \eqref{eql.Carl4_3} becomes
\begin{align}
\label{eql.Carl4_4} \left[ \int_{ \mc{G} } ( \cos 2 \theta )^{ (-1 + 2a) \frac{n}{2} } \right]^\frac{2}{n} &\lesssim \left[ t_\ast^n \int_{- \frac{\pi}{4} }^\frac{\pi}{4} ( \cos 2 \theta )^{ (-1 + 2a) \frac{n}{2} } d \theta \right]^\frac{2}{n} \lesssim t_\ast^2 \text{,}
\end{align}
since \eqref{eql.ass_a} holds, and since $\cos (2 \theta)$ vanishes to first-order at $| \theta | = \frac{\pi}{4}$.

Thus, combining \eqref{eql.Carl4_1}, \eqref{eql.Carl4_2}, and \eqref{eql.Carl4} yields
\begin{align}
\label{eql.Carl4_6} t_\ast^{-1} \int_{ \mc{G} } ( \cos 2 \theta )^{-1 + 2a} |\phi|^2 &\lesssim t_\ast \int_{ \mc{K}^\sigma_\eta (t_\ast) } ( | \nabla \phi |^2 + t_\ast^{-2} | \phi |^2 ) \text{.}
\end{align}
The inequality \eqref{eql.Carl4}, along with \eqref{eql.Carl4_6}, implies the desired estimate \eqref{eq.main_est_timecone}.
\medskip

\noindent
\emph{The case $n \le 2$:}
In this case, we have $H^1 \hookrightarrow L^q$ for any $2 \leq q < \infty$.
Taking $q$ to be sufficiently large so that $\frac{q (-1 + 2a)}{q - 1} > -1$, analogous computations yield
\footnote{Here, we can take any $a$ satisfying the first condition in \eqref{eq.carleman_timecone_ass}.}
\begin{align}
\label{eql.Carl4_7} t_\ast^{-1} \int_{ \mc{G} } ( \cos 2 \theta )^{-1 + 2a} |\phi|^2 &\lesssim t_\ast^{-1} \left[ \int_{ \mc{G} } ( \cos 2 \theta )^\frac{ q (-1 + 2a) }{q - 1} \right]^\frac{q - 1}{q} \left( \int_{ \mc{G} } |\phi|^{2q} \right)^\frac{1}{q} \\
\notag &\lesssim t_\ast^{-1} \cdot t_\ast^{\frac{ n (q - 1) }{q}} \cdot t_\ast^{ \frac{n}{q} - n + 2 } \int_{ \mc{G} } ( | \nabla \phi |^2 + t_\ast^{-2} | \phi |^2 ) \\
\notag &\lesssim t_\ast \int_{ \mc{K}^\sigma_\eta (t_\ast) } ( | \nabla \phi |^2 + t_\ast^{-2} | \phi |^2 ) \text{,}
\end{align}
Combining \eqref{eql.Carl4} with \eqref{eql.Carl4_7} once again results in \eqref{eq.main_est_timecone}.

\subsubsection{Proof of \eqref{eq.main_est_annulus}}

We once again begin from \eqref{eql.Carl4}.
We integrate this equation with respect to $\sigma$ over the interval $[\sigma_0, \sigma_1]$, which yields
\begin{align}
\label{eql.Carl4i_1} \int_{ \mc{R}^\sigma_\gamma (t_\ast) } |\phi|^{p+1} &\lesssim t_\ast \int_{ \sigma_0 }^{ \sigma_1 } \int_{ \mc{K}^\sigma_\eta (t_\ast) } [ ( \partial_t \phi )^2 + ( \partial_r \phi )^2 + | \phi |^{p + 1} + t_\ast^{-2} | \phi |^2 ] d \sigma \\
\notag &\qquad + t_\ast^{-1} \int_{ \sigma_0 }^{ \sigma_1 } \int_{ \mc{G}^\sigma_{t_\ast, \zeta_0} } ( \cos 2 \theta_{t_\ast, \zeta_0} )^{-1 + 2a} |\phi|^2 d \sigma \\
\notag &\qquad + t_\ast^{-1} \int_{ \sigma_0 }^{ \sigma_1 } \int_{ \mc{G}^\sigma_{t_\ast, \zeta_1} } ( \cos 2 \theta_{t_\ast, \zeta_1} )^{-1 + 2a} |\phi|^2 d \sigma \\
\notag &:= I + J_0 + J_1 \text{.}
\end{align}

For the first term on the right-hand side of \eqref{eql.Carl4i_1}, we resort to the coarea formula.
In particular, we think of $\mc{K}^\sigma_\eta (t_\ast)$ as the level set $\frac{r}{t} = \sigma$.
Letting
\begin{align*}
\mc{R}^{ \sigma_0, \sigma_1 }_\eta (t_\ast) := \bigcup_{ \sigma_0 < \sigma < \sigma_1 } \mc{K}^\sigma_\eta (t_\ast) = \mc{R}^{ \sigma_1 }_\eta \setminus \bar{\mc{R}}^{ \sigma_0 }_\eta \text{,}
\end{align*}
and applying a short computation, we see that this can be estimated
\begin{align}
\label{eql.Carl4i_2} I &\lesssim \int_{ \mc{R}^{\sigma_0, \sigma_1}_\eta (t_\ast) } [ ( \partial_t \phi )^2 + ( \partial_r \phi )^2 + | \phi |^{p + 1} + t_\ast^{-2} | \phi |^2 ] \\
\notag &\lesssim t_\ast \sup_{ \eta^{-1} t_\ast \leq \tau \leq \eta t_\ast } \int_{ \mc{A}^{\sigma_0, \sigma_1} (\tau) } [ ( \partial_t \phi )^2 + ( \partial_r \phi )^2 + | \phi |^{p + 1} + t_\ast^{-2} | \phi |^2 ] \text{.}
\end{align}

As the remaining two terms in \eqref{eql.Carl4i_1} will be treated identically, we again fix $i \in \{ 0, 1 \}$, and we abbreviate $\mc{G}^\sigma := \mc{G}^\sigma_{t_\ast, \zeta_i}$ and $\theta := \theta_{t_\ast, \zeta_i}$.
Applying the coarea formula in the same manner as in \eqref{eql.Carl4_3}, we obtain
\begin{align}
\label{eql.Carl4i_3} J_i &\lesssim \int_{ -\frac{\pi}{4} }^\frac{\pi}{4} ( \cos 2 \theta' )^{-1 + 2a} \int_{ \sigma_0 }^{ \sigma_1 } \int_{ (\theta, \sigma) = (\theta', \sigma') } | \phi |^2 d \sigma' d \theta' \\
\notag &\lesssim \sup_{ | \theta' | \leq \frac{\pi}{4} } \int_{ \sigma_0 }^{ \sigma_1 } \int_{ (\theta, \sigma) = (\theta', \sigma') } | \phi |^2 d \sigma' \text{.}
\end{align}
To control the right-hand side of \eqref{eql.Carl4i_3}, we apply the fundamental theorem of calculus by integrating along each $\mc{G}^{\sigma'}$ from $\theta = \theta'$ to $\theta = 0$ (i.e., to $t = t_\ast$):
\begin{align}
\label{eql.Carl4i_4} J_i &\lesssim \int_{ \sigma_0 }^{ \sigma_1 } \int_{ (t, \sigma) = (t_\ast, \sigma') } | \phi |^2 d \sigma' + \int_{ \sigma_0 }^{ \sigma_1 } \int_{ \mc{G}^{\sigma'} } [ | \phi \partial_t \phi | + | \phi \partial_r \phi | + t_\ast^{-1} \phi^2 ] d \sigma' \\
\notag &\lesssim \int_{ \sigma_0 }^{ \sigma_1 } \int_{ (t, \sigma) = (t_\ast, \sigma') } | \phi |^2 d \sigma' + \int_{ \sigma_0 }^{ \sigma_1 } \int_{ \mc{K}^\sigma_\eta (t_\ast) } [ t_\ast | \nabla \phi |^2 + t_\ast^{-1} \phi^2 ] d \sigma' \text{.}
\end{align}

Applying the appropriate coarea formulas to \eqref{eql.Carl4i_4} yields
\begin{align}
\label{eql.Carl4i_5} J_i &\lesssim t_\ast^{-1} \int_{ \mc{A}^{\sigma_0, \sigma_1} (t_\ast) } | \phi |^2 + \int_{ \mc{R}^{\sigma_0, \sigma_1}_\eta (t_\ast) } [ | \nabla \phi |^2 + t_\ast^{-2} \phi^2 ] \\
\notag &\lesssim t_\ast \sup_{ \eta^{-1} t_\ast \leq \tau \leq \eta t_\ast } \int_{ \mc{A}^{\sigma_0, \sigma_1} (\tau) } [ ( \partial_t \phi )^2 + ( \partial_r \phi )^2 + t_\ast^{-2} \phi^2 ] \text{.}
\end{align}
Finally, combining \eqref{eql.Carl4i_1}, \eqref{eql.Carl4i_2}, and \eqref{eql.Carl4i_5} results in the desired estimate \eqref{eq.main_est_annulus}.

\section{Proofs of the Main Theorems} \label{sec:proofs}

In this section, we prove the main results of the paper: Theorems \ref{singularity.foc} and \ref{out-in-decay1}.

\subsection{Proof of Theorem \ref{singularity.foc}.}

Assume the hypotheses of Theorem \ref{singularity.foc}, and let $\delta$ be as in the assumptions.
Furthermore, for convenience, using the time reflection symmetry $t \mapsto -t$, we can work with positive times instead, that is, we can replace $\phi (t, x)$ and $V (t, x)$ by $\phi (-t, x)$ and $V (-t, x)$, respectively.

We begin by applying the estimate from Theorem \ref{thm.main_est_annulus}.
First, we observe that there is some $0 < t_0 \ll 1$ such that for any $0 < t_\ast < t_0$:
\begin{itemize}
\item Applying \eqref{eq.singularity_ass} and the $H^1$-$L^{p+1}$-Sobolev inequality, we have
\begin{align}
\label{eql.singularity_0a} t_\ast^{2 - n + \frac{4}{p-1}} \int_{ \mc{A}^{\sigma_0, \sigma_1} (t_\ast) } ( | \nabla \phi |^2 + | \phi |^{p + 1} + t_\ast^{-2} |\phi|^2 ) \lesssim \delta \text{.}
\end{align}

\item The conditions for $V$ in \eqref{eq.singularity_ass_pV} imply that \eqref{eq.main_est_ass_V} holds for $t_\ast$.
\end{itemize}
Thus, \eqref{eq.main_est_annulus} yields, for $t_\ast < t_0$ and for some $\eta > \gamma > 1$, that
\begin{align}
\label{eql.singularity_0} \int_{ \mc{R}^{\sigma_0}_\gamma (t_*) } |\phi|^{p+1} &\lesssim t_\ast \sup_{ \eta^{-1} t_\ast \leq \tau \leq \eta t_\ast } \int_{ \mc{A}^{\sigma_0, \sigma_1}_\eta (\tau) } ( | \nabla \phi |^2 + |\phi|^{p+1} + t_\ast^{-2} | \phi |^2 ) \\
\notag &\lesssim \delta t_\ast^{n - 1 - \frac{4}{p - 1}} \text{.}
\end{align}
Note that $n - 1 - \frac{4}{p - 1} < 0$, since $p$ is subconformal.

Moreover, combining \eqref{eql.singularity_0} with H\"older's inequality yields
\begin{align}
\label{eql.singularity_1} \int_{ \mc{R}^{\sigma_0}_\gamma (t_\ast) } \phi^2 \lesssim \delta t_\ast^{n + 1 - \frac{4}{p-1}} \text{,} \qquad 0 < t_\ast < t_0 \text{.}
\end{align}
Thus, it remains only to show that
\begin{align}
\label{eql.singularity_2} \limsup_{ t_\ast \searrow 0 } t_\ast^{n - 1 - \frac{4}{p-1}} \int_{ \mc{R}^{\sigma_0}_\gamma (t_\ast) } | \nabla \phi |^2 \lesssim \delta \text{.}
\end{align}
If \eqref{eql.singularity_2} holds, then \eqref{eql.singularity_1} and \eqref{eql.singularity_2} imply the desired conclusion, \eqref{eq.singularity}.

\subsubsection{Proof of \eqref{eql.singularity_2}}

Fix $0 < t_\ast < t_0$, and consider the domain
\begin{align*}
\Theta &:= \bigcup_{ t_\ast < \tau < 1 } \mc{A}^{\sigma_0, \sigma_1} (\tau) = \bigcup_{ \sigma_0 < \sigma < \sigma_1 } ( \partial \mc{C}^\sigma \cap \{ t_\ast < t < 1 \} ) \text{.}
\end{align*}
Applying the coarea formula (similar to \eqref{eql.Carl4i_1}), we see that
\begin{align}
\label{eql.singularity_30} &\int_{ \sigma_0 }^{ \sigma_1 } \int_{ \partial \mc{C}^\sigma \cap \{ t_\ast < t < 1 \} } t ( | \nabla \phi |^2 + | \phi |^{p + 1} + t^{-2} \phi^2 ) d \sigma \\
\notag &\qquad \lesssim \int_{t_\ast}^1 \int_{ \mc{A}^{\sigma_0, \sigma_1} (t) } ( | \nabla \phi |^2 + | \phi |^{p + 1} + t^{-2} \phi^2 ) dt \\
\notag &\qquad \lesssim C_\delta + \delta t_\ast^{n - 1 - \frac{4}{p - 1}} \text{,}
\end{align}
where $C_\delta$ denotes a constant which depends on $\delta$.
In particular, the portion of the integral in $t_\ast < t < \delta$ is controlled using \eqref{eql.singularity_0a}, while the part $t_0 < t < 1$ is trivially bounded by a constant.
It follows that there is some $\sigma' \in (\sigma_0, \sigma_1)$ such that
\begin{align}
\label{eql.singularity_3} t_\ast \int_{ \partial \mc{C}^{\sigma'} \cap \{ t_\ast < t < 1 \} } ( | \nabla \phi |^2 + | \phi |^{p + 1} + t^{-2} \phi^2 ) &\lesssim C_\delta + \delta t_\ast^{n - 1 - \frac{4}{p - 1}} \text{.}
\end{align}

Fix $t_1 < \eta^{-1}$, and consider now the region $\Omega_{t_\ast}$ bounded by the hypersurfaces:
\begin{align*}
\{ t = t_1 \} \text{,} \qquad \{ t = t_\ast \} \text{,} \qquad \partial \mc{C}^{\sigma'} \cap \{ t_\ast < t < 1 \} \text{.}
\end{align*}
By similar reasoning as before, and by \eqref{eql.singularity_0a}, we can bound
\begin{align}
\label{eql.singularity_40} \int_{ \Omega_{t_\ast} } | \phi |^{p + 1} &\lesssim C_\delta + \int_{ \Omega_{t_\ast} \cap \{ t < t_0 \} } | \phi |^{p + 1} \\
\notag &\lesssim C_\delta + \int_{t_\ast}^{t_0} \int_{ \mc{A}^{\sigma_0, \sigma_1} (t) } | \phi |^{p + 1} dt + \int_{ \mc{C}^{\sigma_0} \cap \{ t_\ast < t < t_0 \} } | \phi |^{p + 1} \\
\notag &\lesssim C_\delta + \delta t_\ast^{n - 1 - \frac{4}{p - 1}} + \int_{ \mc{C}^{\sigma_0} \cap \{ t_\ast < t < t_0 \} } | \phi |^{p + 1} \text{.}
\end{align}
The integral on the right-hand side of \eqref{eql.singularity_40} can be decomposed into slabs of the form $\mc{R}^{\sigma_0}_\gamma (\gamma^k t_\ast)$, which can be controlled using \eqref{eql.singularity_0}.
Summing these, we see that
\begin{align}
\label{eql.singularity_4} \int_{ \Omega_{t_\ast} } | \phi |^{p + 1} &\lesssim C_\delta + \delta t_\ast^{n - 1 - \frac{4}{p - 1}} \text{.}
\end{align}

We now resort to the standard energy identity for the wave operator.
Applying the divergence theorem in the usual way on $\Omega_{t_\ast}$, we obtain
\begin{align}
\label{eql.singularity_50} \int_{ \Omega_{t_\ast} } \Box \phi \partial_t \phi &= \frac{1}{2} \int_{ |x| \leq \sigma' t_\ast } | \nabla \phi (t_\ast, x) |^2 dx - \frac{1}{2} \int_{ |x| \leq \sigma' t_1 } | \nabla \phi (t_1, x) |^2 dx \\
\notag &\qquad + \int_{ \partial \mc{C}^{\sigma'} \cap \{ t_\ast < t < 1 \} } K_1 \text{,}
\end{align}
where the function $K_1$ satisfies
\begin{align}
\label{eql.singularity_51} | K_1 | \lesssim | \nabla \phi |^2 \text{.}
\end{align}

Next, since
\begin{align}
\label{eql.singularity_52} \int_{ \Omega_{t_\ast} } \Box \phi \partial_t \phi &= - \frac{1}{p + 1} \int_{ \Omega_{t_\ast} } V \partial_t | \phi |^{p + 1} \text{,}
\end{align}
then applying the divergence theorem to the right-hand side of \eqref{eql.singularity_52}, recalling \eqref{eql.singularity_0} and \eqref{eql.singularity_1}, and recalling that $V$ and $\partial_t V$ are uniformly bounded, we obtain
\begin{align}
\label{eql.singularity_53} \int_{ |x| \leq \sigma' t_\ast } | \nabla \phi (t_\ast) |^2 &\lesssim \int_{ |x| \leq \sigma' t_\ast } | \phi (t_\ast) |^{p+1} + \int_{ |x| \leq \sigma' t_1 } [ | \nabla \phi (t_1) |^2 + | \phi (t_1) |^{p + 1} ] \\
\notag &\qquad + \int_{ \partial \mc{C}^{\sigma'} \cap \{ t_\ast < t < 1 \} } ( | \nabla \phi |^2 + | \phi |^{p + 1} ) + \int_{\Omega_{t_\ast}} | \phi |^{p + 1} \text{.}
\end{align}
Recalling \eqref{eql.singularity_3} and \eqref{eql.singularity_4} yields
\begin{align}
\label{eql.singularity_5} \int_{ |x| \leq \sigma' t_\ast } | \nabla \phi (t_\ast) |^2 &\lesssim \int_{ |x| \leq \sigma' t_\ast } | \phi (t_\ast) |^{p+1} + C_\delta t_\ast^{-1} + \delta t_\ast^{n - 2 - \frac{4}{p - 1}} \text{.}
\end{align}

Applying \eqref{eql.singularity_5} for all times $\gamma^{-1} t_\ast < t < \gamma t_\ast$ and integrating over $t$, we have
\begin{align}
\label{eql.singularity_6} \int_{ \mc{R}^{ \sigma_0 }_\gamma (t_\ast) } | \nabla \phi (t_\ast) |^2 &\lesssim \int_{ \mc{R}^{ \sigma_0 }_\gamma (t_\ast) } | \phi (t_\ast) |^{p+1} + C_\delta + \delta t_\ast^{n - 1 - \frac{4}{p - 1}} \\
\notag &\lesssim C_\delta + \delta t_\ast^{n - 1 - \frac{4}{p - 1}} \text{.}
\end{align}
Finally, since $n - 1 - \frac{4}{p - 1} < 0$, then letting $t_\ast \searrow 0$ in \eqref{eql.singularity_6} results in \eqref{eql.singularity_2}.
 

\subsection{Proof of Theorem \ref{out-in-decay1}}

First, H\"older's inequality implies
\begin{align*}
\int_{ \partial \mc{C}^\sigma \cap \{ t > 1 \} } t^{-2} \phi^2 &\lesssim \left( \int_{ \partial \mc{C}^\sigma \cap \{ t > 1 \} } |\phi|^{p+1} \right)^\frac{2}{p + 1} \left( \int_{ \partial \mc{C}^\sigma \cap \{ t > 1 \} } t^{ -2 \frac{p+1}{p-1} } \right)^\frac{p-1}{p+1} \\
&\lesssim \left( \int_{ \partial \mc{C}^\sigma \cap \{ t > 1 \} } |\phi|^{p+1} \right)^\frac{2}{p + 1} \left( \int_1^\infty t^{-2 - \frac{4}{p-1} + n - 1} dt \right)^\frac{p - 1}{p + 1} \text{.}
\end{align*}
Since $p < 1 + \frac{4}{n - 1}$, the last integral in the above is finite, thus
\begin{align}
\label{Holder-appl} \int_{ \partial \mc{C}^\sigma \cap \{ t > 1 \} } ( | \nabla \phi |^2 + | \phi |^{p + 1} + t^{-2} \phi^2 ) < \infty.
\end{align}

By the assumptions \eqref{eq.final_state_ass_V} for $V$, there exists $t_0 > 0$ such that \eqref{eq.main_est_ass_V} holds for any $t_\ast > t_0$, with sufficiently small $\alpha$.
Fixing now $t_1 \gg t_0$, it suffices to show that
\begin{align}
\label{eql.final_state_0} \int_{ \mc{C}^\sigma \cap \{ t > t_1 \} } t^{-1} | \phi |^{p + 1} < \infty \text{,}
\end{align}
since the corresponding integral on $\{ 1 < t \leq t_1 \}$ is trivially finite.

The key step is to apply Theorem \ref{thm.main_est_timecone}: there exist $\eta > \gamma > 1$, as in the statement of Theorem \ref{thm.main_est_timecone},\ such that we have from \eqref{eq.main_est_timecone} the estimate
\begin{align}
\label{eql.final_state_1} ( \gamma^k t_1 )^{-1} \int_{ \mc{R}^\sigma_\gamma (\gamma^k t_1) } |\phi|^{p+1} \lesssim \int_{ \mc{K}^\sigma_\eta (\gamma^k t_1) } [ | \nabla \phi |^2 + |\phi|^{p+1} + ( \gamma^k t_1 )^{-2} | \phi |^2 ] \text{,}
\end{align}
for any $k > 0$.
We then cover $\mc{C}^\sigma \cap \{ t > t_1 \}$ by a tiling of slabs of the form $\mc{R}^\sigma_\gamma (\gamma^k t_1)$, for all $k > 0$.
Observe there is a number $M > 0$, depending on $\eta$ and $\gamma$, so that the each point of $\partial \mc{C}^\sigma \cap \{ t > 1 \}$ intersects at most $M$ sets of the form $\mc{K}^\sigma_\eta (\gamma^k t_1)$.
Thus, applying \eqref{eql.final_state_1} for $k > 0$ and summing over the $k$, we obtain
\begin{align}
\label{eql.final_state_2} \int_{ \mc{C}^\sigma \cap \{ t > t_1 \} } t^{-1} |\phi|^{p+1} \lesssim \int_{ \partial \mc{C}^\sigma \cap \{ t > 1 \} } ( | \nabla \phi |^2 + | \phi |^{p+1} + t^{-2} | \phi |^2 ) < \infty \text{.}
\end{align}
This establishes \eqref{eql.final_state_0} and completes the proof of Theorem \ref{out-in-decay1}.

\raggedright
\bibliographystyle{amsplain}
\bibliography{bib}

\end{document}